\documentclass[11pt]{article}
\usepackage{cite}
\usepackage{mathrsfs}
\usepackage{amsfonts}
\usepackage{amsmath}
\usepackage{amsfonts,amssymb,color}
\usepackage{dsfont}
\usepackage{curves}
\usepackage{mathrsfs}
\usepackage{pifont}
\usepackage{amssymb}
\usepackage{latexsym,amsmath,amssymb,amsfonts,epsfig,graphicx,cite,psfrag}
\usepackage{eepic,color,colordvi,amscd}

\newtheorem{theorem}{Theorem}[section]
\newtheorem{proposition}{Proposition}[section]

\newtheorem{lemma}{Lemma}[section]

\newtheorem{problem}{Problem}[section]
\newtheorem{claim}{Claim}[section]
\newtheorem{conjecture}{Conjecture}[section]

\newcommand{\qed}{\hfill\rule{0.5em}{0.809em}}

\def\emptyset{\mbox{{\rm \O}}}
\def\bar{\overline}

\def\qed{\hfill \rule{4pt}{7pt}}

\def\pf{\noindent {\it Proof. }}

\begin{document}
\title{On the chromatic number of some ($P_3\cup P_2$)-free graphs\thanks{Partially supported by Natural Science Foundation of Jiangsu Province (No. BK20170862) and NSFC 11701142.}}

\author{{Rui Li $^{a,*}$\footnote{Email address: lirui@hhu.edu.cn}\quad Jinfeng Li $^{a,}$\footnote{Email address: 1345770246@qq.com}\quad  Di Wu$^{b,}$\footnote{Email address: 1975335772@qq.com
}}\\
       {\small $^a$ School of Mathematics, Hohai University}\\
        {\small 8 West Focheng Road, Nanjing, 211100, China}\\
        {\small $^b$ School of Mathematical Science, Nanjing Normal University}\\
        {\small 1 Wenyuan Road, Nanjing, 210046, China}
        }
\date{}
\maketitle

\begin{abstract}
A hereditary class $\cal G$ of graphs is {\em $\chi$-bounded} if there is a {\em $\chi$-binding function}, say $f$, such that $\chi(G)\le f(\omega(G))$ for every $G\in\cal G$, where $\chi(G)(\omega(G))$ denotes the chromatic (clique) number of $G$. It is known that for every $(P_3\cup P_2)$-free graph $G$, $\chi(G)\le \frac{1}{6}\omega(G)(\omega(G)+1)(\omega(G)+2)$ \cite{BA18}, and the class of $(2K_2, 3K_1)$-free graphs does not admit a linear $\chi$-binding function\cite{BBS19}. In this paper, we prove that (\romannumeral 1) $\chi(G)\le2\omega(G)$ if $G$ is ($P_3\cup P_2$, kite)-free, (\romannumeral 2) $\chi(G)\le\omega^2(G)$ if $G$ is ($P_3\cup P_2$, hammer)-free, (\romannumeral 3) $\chi(G)\le\frac{3\omega^2(G)+\omega(G)}{2}$ if $G$ is ($P_3\cup P_2, C_5$)-free. Furthermore, we also discuss $\chi$-binding functions for $(P_3\cup P_2, K_4)$-free graphs.

\begin{flushleft} {\em Key words and
phrases:}  chromatic number; clique number; $\chi$-binding function; $(P_3\cup P_2)$-free graphs 

\end{flushleft}
\end{abstract}

\section{Introduction}

All graphs considered in this paper are finite and simple. We use $P_k$ and $C_k$ to denote a {\em path} and a {\em cycle} on $k$ vertices respectively, and follow \cite{BM1976} for undefined notations and terminology. Let $G$ be a graph, and $X$ be a subset of $V(G)$. We use $G[X]$ to denote the subgraph of $G$ induced by $X$, and call $X$ a {\em clique} ({\em independent set}) if $G[X]$ is a complete graph (has no edge). The {\em clique number} $\omega(G)$ of $G$ is the maximum size taken over all cliques of $G$.

For $v\in V(G)$, let $N_G(v)$ be the set of vertices adjacent to $v$, $d_G(v)=|N_G(v)|$, $N_G[v]=N_G(v)\cup \{v\}$, $M_G(v)=V(G)\setminus N_G[v]$. For $X\subseteq V(G)$, let $N_G(X)=\{u\in V(G)\setminus X\;|\; u$ has a neighbor in $X\}$ and $M_G(X)=V(G)\setminus (X\cup N_G(X))$. If it does not cause any confusion, we will omit the subscript $G$ and simply write $N(v), d(v), N[v], M(v), N(X)$ and $M(X)$.  Let $\delta(G)$ denote the minimum degree of $G$.

For positive integer $i$, let $N^i_G(X):=\{u\in V(G)\setminus X\ | \ \min \{d_G(u,v)\}=i\ for\ v\in X$\}, where $d_G(u,v)$ is the distance between $u$ and $v$ in $G$. Then $N^1_G(X)=N_G(X)$ is the neighborhood of $X$ in $G$. Moreover, let $N^{\ge i}_G(X):= \cup^\infty_{j=i}N^i_G(X)$. We write $N^i_G(H)$ for $N^i_G(V(H))$.

Let $G$ and $H$ be two vertex disjoint graphs. The {\em union} $G\cup H$ is the graph with $V(G\cup H)=V(G)\cup (H)$ and $E(G\cup H)=E(G)\cup E(H)$. The union of $k$ copies of the same graph $G$ will be denoted by $kG$. The {\em join} $G+H$ is the graph with $V(G+H)=V(G)\cup V(H)$ and $E(G+H)=E(G)\cup E(H)\cup\{xy\;|\; x\in V(G), y\in V(H)$$\}$. The complement of a graph $G$ will be denoted by $\bar G$.  We say that $G$ induces $H$ if $G$ has an induced subgraph isomorphic to $H$, and say that $G$ is $H$-free otherwise. Analogously, for a family $\cal H$ of graphs, we say that $G$ is ${\cal H}$-free if $G$ induces no member of ${\cal H}$.

Let $k$ be a positive integer, and let $[k]=\{1, 2, \ldots, k\}$. A $k$-{\em coloring} of $G$ is a mapping $c: V(G)\mapsto [k]$ such that $c(u)\neq c(v)$ whenever $u\sim v$ in $G$. The {\em chromatic number} $\chi(G)$ of $G$ is the minimum integer $k$ such that $G$ admits a $k$-coloring.  It is certain that $\chi(G)\ge \omega(G)$. A {\em perfect graph} is one such that  $\chi(H)=\omega(H)$ for all of its induced subgraphs $H$.  A family $\cal G$ of graphs is said to be $\chi$-{\em bounded} if there is a function $f$ such that $\chi(G)\le f(\omega(G))$ for every $G\in\cal G$, and if such a function does exist for $\cal G$, then $f$ is said to be a {\em binding function} of $\cal G$ \cite{Gy75}.

Let $G$ be a graph, $v\in V(G)$, and let $X$ and $Y$ be two subsets of $V(G)$. We say that $v$ is {\em complete} to $X$ if $v$ is adjacent to all vertices of $X$, and say that $v$ is {\em anticomplete} to $X$ if $v$ is not adjacent to any vertex of $X$. We say that $X$ is complete (resp. anticomplete) to $Y$ if each vertex of $X$ is complete (resp. anticomplete) to $Y$.  Particularly, we say that $X$ is {\em almost complete}  to $Y$ if at most one vertex of $X$ is not complete to $Y$. For $u, v\in V(G)$, we simply write $u\sim v$ if $uv\in E(G)$, and write $u\not\sim v$ if $uv\not\in E(G)$.

A {\em hole} of $G$ is an induced cycle of length at least 4, and a {\em $k$-hole} is a hole of length $k$. A $k$-hole is called an {\em odd hole} if $k$ is odd, and is called an {\em even hole} otherwise. An {\em antihole} is the complement of some hole. An odd (resp. even) antihole is defined analogously. The famous {\em Strong Perfect Graph Theorem} states that

\begin{theorem}\label{Perfect}\cite{CRSR06}
	A graph is perfect if and only if it induces neither an odd hole nor an odd antihole.
\end{theorem}

Erd\"{o}s \cite{ER59} proved that for any positive integers $k,l\ge3$, there exists a graph $G$ with $\chi(G)\ge k$ and no cycles of length less than $l$. This result motivates us to study the chromatic number of $F$-free graphs, where $F$ is a forest (a disjoint union of trees). Gy\'{a}rf\'{a}s \cite{Gy75} and Sumner \cite{Su81} independently, proposed the following famous conjecture.

\begin{conjecture}\label{tree}\cite{Gy75, Su81}
	Let $F$ be a forest. Then $F$-free graphs are $\chi$-bounded.
\end{conjecture}

The class of $2K_2$-free graphs has attracted a great deal of interest in recent years.  It is known that for every $2K_2$-free graph $G$, $\chi(G)\le \binom{\omega+1}{2}$ \cite{Wa80}. Up to now, the best known $\chi$-binding function for $2K_2$-free graphs is $f(\omega)=\binom{\omega+1}{2}-2\lfloor \frac{\omega}{3} \rfloor$ \cite{GM22}. We refer the interested readers to \cite{BBS19,KM18,PA23,Wa80}  for results of $2K_2$-free graphs, and to \cite{RS04,SR19,SS20} for more results and problems about the $\chi$-bounded problem. In particular, Brause {\em et al.}\cite{BBS19} proved that the class of $(2K_2, 3K_1)$-free graphs does not admit a linear $\chi$-binding function. In 2022, Brause {\em et al.}\cite{BMS19} gave a more general theorem.

\begin{theorem}\label{good}(Lemma 1 of \cite{BM1976})
	Let $\cal H$ be a set of graphs and $\ell$ be an integer such that $\bar{H}$ has girth at most $\ell$ for each $H\in\cal H$. If the class of $\cal H$-free graphs is $\chi$-bounded, then the class of $\cal H$-free graphs does not admit a linear $\chi$-binding function.
\end{theorem}

\begin{figure}[htbp]\label{fig-1}
	\begin{center}
		\includegraphics[width=15cm]{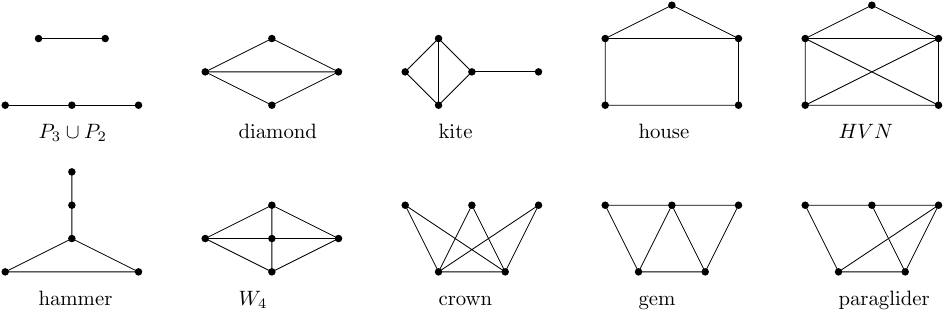}
	\end{center}
	\vskip -15pt
	\caption{Illustration of $P_3\cup P_2$ and some forbidden configurations.}
\end{figure}

Since $(P_3\cup P_2)$-free graphs is a superclass of $2K_2$-free graphs, many scholars began to show interest in $(P_3\cup P_2)$-free graphs. The best known $\chi$-binding function for $(P_3\cup P_2)$-free graphs is $f(\omega)=\frac{1}{6}\omega(\omega+1)(\omega+2)$ \cite{BA18}. In \cite{PA22}, Prashant {\em  et al.} proved that if $G$ is ($P_3\cup P_2$, diamond)-free, then $\chi(G)\le4$ when $\omega(G)=2$, $\chi(G)\le6$ when $\omega(G)=3$, $\chi(G)\le4$ when $\omega(G)=4$, and $G$ is perfect when $\omega(G)\ge5$, and they also proved \cite{PA22} that $\chi(G)\le\omega(G)+1$ if $G$ is a $(P_3\cup P_2, HVN)$-free graph with $\omega(G)\ge4$.  As a superclass of ($P_3\cup P_2$, diamond)-free graphs, Cameron {\em  et al.} \cite{CHM21} proved that $\chi(G)\le\omega(G)+3$ if $G$ is ($P_6$, diamond)-free, this bound is optimal.  Very Recently, Wu and Xu \cite{WX23} proved that $\chi(G)\le\frac{1}{2}\omega^2(G)+\frac{3}{2}\omega(G)+1$ if $G$ is ($P_3\cup P_2$, crown)-free, Char and Karthick \cite{CK22} proved that $\chi(G)\le$ max $\{\omega(G)+3, \lfloor \frac{3\omega(G)}{2} \rfloor-1\}$ if $G$ is a ($P_3\cup P_2$, paraglider)-free graph with $\omega(G)\ge3$, Prashant {\em et al.} proved that $\chi(G)\le2\omega(G)$ if $G$ is ($P_3\cup P_2$, gem)-free, and  Li {\em et al.}\cite{LLW23,LLWD23} proved that $\chi(G)\le2\omega(G)$ if $G$ is ($P_3\cup P_2$, house)-free or ($P_3\cup P_2, W_4$)-free. (See Figure \ref{fig-1} for the illustration of $P_3\cup P_2$ and some forbidden configurations.)

In this paper, we prove that

\begin{theorem}\label{kite}
	$\chi(G)\le2\omega(G)$ if $G$ is $($$P_3\cup P_2$, kite$)$-free.
\end{theorem}

Let $G$ be a graph on $n$-vertices $\{v_1, v_2, \dots , v_n\}$ and let $H_1, H_2, \dots , H_n$ be $n$ vertex-disjoint graphs. An expansion
$G(H_1, H_2, \dots , H_n)$ of $G$ is a graph obtained from $G$ by (i) replacing each $v_i$ of $G$ by $H_i$, $i = 1, 2, . . . , n$ and
(ii) by joining every vertex in $H_i$ with every vertex in $H_j$, whenever $v_i$ and $v_j$ are adjacent in $G$. In addition, $G\cong K_n$ and $H_i\cong H$, denote by $K_n(H)=G(H,H,\dots, H)$.

Let $H$ be the Mycielski-Gr\"{o}stzsch graph (see Figure 2). Then $\omega(H)=2$ and $\chi(H)=4$. It is clear that $K_k(H)$ is ($P_3\cup P_2$, kite)-free, and $\chi(K_k(H))=2\omega(K_k(H))=2\cdot 2k=4k$. Let $H'$ be the complement of Schl\"{a}fli graph (see https://houseofgraphs.org/graphs/19273). Then $\omega(H')=3$ and $\chi(H')=6$. Obviously, $K_{k-1}(H)+H'$ is ($P_3\cup P_2$, kite)-free. Hence $\omega(K_{k-1}(H)+H')=2k+1$ and $\chi(K_{k-1}(H)+H')=4k+2$. This implies that our bound is optimal when $\omega(G)\ge2$.

By Theorem \ref{good}, we have that the class of ($P_3\cup P_2$, hammer)-free graphs does not admit a linear $\chi$-binding function. In this paper, we prove that

\begin{theorem}\label{hammer}
	$\chi(G)\le\omega^2(G)$ if $G$ is $($$P_3\cup P_2$, hammer$)$-free.
\end{theorem}

Notice that $H$ is also ($P_3\cup P_2$, hammer)-free. Therefore, the $\chi$-binding function $f=\omega^2$ for ($P_3\cup P_2$, hammer)-free graphs is tight when $\omega=2$.

\begin{figure}[htbp]\label{fig-2}
	\begin{center}
		\includegraphics[width=4cm]{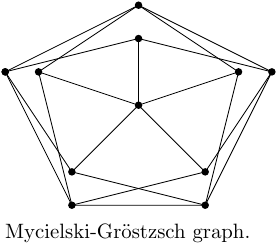}
	\end{center}
	\vskip -15pt
	\caption{Mycielski-Gr\"{o}stzsch graph.}
\end{figure}

In \cite{CK10},  Choudum and Karthick proved that $\chi(G)\le \lceil \frac{5\omega(G)}{4} \rceil$ if $G$ is $(P_3\cup P_2, C_4)$-free. By Theorem \ref{good}, we have that the class of $(P_3\cup P_2, C_5)$-free graphs does not admit a linear $\chi$-binding function. In this paper, we prove that

\begin{theorem}\label{C5}
	$\chi(G)\le\frac{3\omega^2(G)+\omega(G)}{2}$ if $G$ is $($$P_3\cup P_2, C_5$$)$-free. 
\end{theorem}

In \cite{WZ22}, Wang and Zhang proved that if $G$ is a $(P_3\cup P_2, K_3)$-free
graph, then $\chi(G)\le3$ unless $G$ is one of eight graphs with $\Delta(G)=5$ and $\chi(G)=4$. In this paper, we also prove that

\begin{theorem}\label{K4}
	$\chi(G)\le9$ if $G$ is $($$P_3\cup P_2, K_4$$)$-free. 
\end{theorem}

Notice that $H'$ is $(P_3\cup P_2)$-free, and $\omega(H')=3,\chi(H')=6$. Therefore, $\chi(G)\ge6$ if $G$ is $($$P_3\cup P_2, K_4$$)$-free. In reality, as mentioned above, $\chi(G)\le6$ if $G$ is $(P_3\cup P_2, H, K_4)$-free, where $H\in\{$ diamond, paraglider, gem, house, $W_4$, kite$\}$. So, we may ask a question that 

\begin{problem}\label{=6}
	Is that true $\chi(G)\le6$ if $G$ is $($$P_3\cup P_2, K_4$$)$-free?
\end{problem}

We will prove Theorem~\ref{kite} in Section 2, prove Theorem~\ref{hammer} in Section 3, prove Theorem~\ref{C5} in Section 4, and prove Theorem~\ref{K4} in Section 5.

\section{($P_3\cup P_2$, kite)-free graphs}
In this section, we consider ($P_3\cup P_2$, kite)-free graphs. We may always assume that $G$ is a ($P_3\cup P_2$, kite)-free graph such that $\chi(G')\le 2\omega(G')$ for every induced subgraph $G'$ of $G$ different from $G$, and $\chi(G)>2\omega(G)$. The following lemmas will be used in our proof.

\begin{lemma}\label{nonadjcent}(Lemma 2.1 of \cite{LLW23})
	Let $u$ and $v$ be two nonadjacent vertices in $G$. Then $N(u)\not\subseteq N(v)$ and $N(v)\not\subseteq N(u)$.
\end{lemma}
\pf Suppose to its contrary that $N(u)\subseteq N(v)$ by symmetry. By assumption, $\chi(G-u)\le 2\omega(G-u) $. Since we can color $u$ by the color of $v$, it follows that $\chi(G)\le 2\omega(G-u) \le 2\omega(G)$, a contradiction.
\qed

\begin{lemma}\label{C3}\cite{WZ22}
	If $G$ is a $(P_3\cup P_2, C_3)$-free  graph, then $\chi(G)\le 4$.
\end{lemma}

\begin{lemma}\label{K3}(1.8 of \cite{ CS10})
	Let $G$ be a $($$K_1\cup K_3$$)$-free graph. If $G$ contains a $K_3$, then $\chi(G)\le2\omega(G)$.
\end{lemma}

We will complete the proof of Theorem \ref{kite} by the following three claims. 

\begin{claim}\label{c1}
	$G$ is $(P_2\cup K_3)$-free.
\end{claim}
\pf Suppose not. Let $Q$ be an induced $P_2\cup K_3$ in $G$ with $V(Q)=\{u_1,u_2,v_1,v_2,v_3\}$ such that $v_1v_2v_3v_1$ is a triangle. We will prove that

\begin{equation}\label{eqa-1}
	\mbox{$N(u_1)\setminus N(u_2)$ and $N(u_2)\setminus N(u_1)$ are both independent}.
\end{equation}
Suppose $N(u_1)\setminus N(u_2)$ has two adjacent vertices $x_1$ and $x_2$. If $|\{v_1,v_2,v_3\}\setminus N(x_1)|\ge2$, we may by symmetry assume that $x_1\not\sim v_1$ and $x_1\not\sim v_2$, then $\{v_1,v_2,u_1,u_2,x_1\}$ induces a $P_3\cup P_2$, a contradiction. So, $|N(x_1)\cap\{v_1,v_2,v_3\}|\ge2$. By symmetry, $|N(x_2)\cap\{v_1,v_2,v_3\}|\ge2$. Therefore, there must exist a vertex in $\{v_1,v_2,v_3\}$ which is complete to $\{x_1,x_2\}$, say $v_1$. Now, $\{u_1,u_2,x_1,x_2,v_1\}$ induces a kite, a contradiction. So, $N(u_1)\setminus N(u_2)$ is independent, and by symmetry, $N(u_2)\setminus N(u_1)$ is independent. This proves (\ref{eqa-1}).

Let $M=G[M(\{u_1,u_2\})]$. Since $G$ is $(P_3\cup P_2)$-free, we have that $G[M]$ is $P_3$-free. Choose the maximum clique in $G[M]$, say $C_1$. Without loss of generality, we may assume that $\{v_1,v_2,v_3\}\subseteq C_1$. Let $D=\{v\in N(u_1)\cap N(u_2)| v$ is complete to $ N(u_1)\cap N(u_2)\setminus \{v\} \}$, and $C=N(u_1)\cap N(u_2)\setminus D$. So, $D$ is a clique, and for each vertex $y$ in $C$, $y$ is not complete to $C\setminus\{y\}$. 

Suppose there exists a vertex $y_1\in C$ such that $y$ is anticomplete to $C_1$. Since $y_1$ is not complete to $C\setminus\{y_1\}$, it follows that $C$ has a vertex $y_2$ such that $y_2\not\sim y_1$. If $y_2$ is anticomplete to $C_1$, then $\{y_1,u_1,y_2,v_1,v_2\}$ induces a $P_3\cup P_2$, a contradiction. So, $y_2$ has a neighbor in $C_1$, say $v_1$. Then $\{y_1,u_1,u_2,y_2,v_1\}$ induces a kite, a contradiction. So, 

\begin{equation}\label{eqa-2}
	\mbox{$C\subseteq N(C_1)$}.
\end{equation}

Next, we prove that 

\begin{equation}\label{eqa-3}
	\mbox{$C$ is complete to $C_1$}.
\end{equation}

Suppose to its contrary. We may assume that there exists a vertex $y_3$ in $C$ such that $y_3\not\sim v_1$. By (\ref{eqa-2}), $y_3\in N(C_1)$. If $|N(y_3)\cap C_1|\ge2$, let $y_3\sim v_2$ and $y_3\sim v_3$, then $\{u_1,y_3,v_1,v_2,v_3\}$ induces a kite, a contradiction. So, $|N(y_3)\cap C_1|=1$, by symmetry, $y_3\sim v_2$.
Similarly, for each vertex $y\in C$, if $y$ is not complete to $C_1$, then $|N(y)\cap C_1|=1$.

By the definition of $C$, we have that there exists a vertex $y_4$ in $C$ such that $y_3\not\sim y_4$. Suppose $y_4$ is not complete to $C_1$. Then $|N(y_4)\cap C_1|=1$. If $y_4\sim v_2$, then $\{y_3,u_1,y_4,v_1,v_3\}$ induces a $P_3\cup P_2$, a contradiction. If $y_4\not\sim v_2$, then $\{u_1,u_2,y_3,y_4,v_2\}$ induces a kite, a contradiction. So, $y_4$ is complete to $C_1$. Now, $\{u_1,u_2,y_3,y_4,v_1\}$ induces a kite, a contradiction. This proves (\ref{eqa-3}).

Let $\omega_1=\omega(G[C])$. By (\ref{eqa-3}) and $C_1$ is a maximum clique in $G[M]$, we have that $\omega(G[M])\le \omega(G)-\omega_1$. So, $\chi(G[M\cup\{u_1,u_2\}])\le\omega(G)-\omega_1$. It is clear that $D$ is complete to $C$. Since $\omega(G[N(u_1)\cap N(u_2)])\le\omega(G)-2$, we have that $\omega(G[D])\le\omega(G)-2-\omega_1$. 

Note that $V(G)=\{u_1,u_2\}\cup N(\{u_1,u_2\})\cup M$ and $N(\{u_1,u_2\})=(N(u_1)\setminus N(u_2))\cup (N(u_2)\setminus N(u_1))\cup (N(u_1)\cap N(u_2))$. By (\ref{eqa-1}), $\chi(G)\le \chi(G[C])+\chi(G[D])+2+\chi(G[M\cup\{u_1,u_2\}])\le f(\omega_1)+2+2\omega(G)-2\omega_1-2\le2\omega(G)2\omega_1-2\omega_1=2\omega(G)$, a contradiction.

This proves Claim \ref{c1}.\qed

\begin{claim}\label{c2}
	$G$ is hammer-free.
\end{claim}
\pf Suppose not. Let $Q$ be an induced hammer in $G$ with $V(Q)=\{v_1,v_2,v_3,v_4,v_5\}$ such that $v_1v_2v_3v_1$ is a triangle, $v_4\sim v_3$ and $d_Q(v_5)=1$. For a subset $S\subseteq\{1,2,4,5\}$, we define $N_S=\{v|v\in N(Q\setminus\{v_3\}),$ and $v_i\sim v$ if and only if $i\in S\}$. Let $N=N(\{v_1,v_2,v_4,v_5\})$ and $M=M(\{v_1,v_2\})$. By Claim \ref{c1}, we have that $G$ is $(P_2\cup K_3)$-free.

If $N_{\{1,2\}}\ne\emptyset$, let $v_1'\in N_{\{1,2\}}$, then $\{v_1,v_2,v_1',v_4,v_5\}$ induces a $P_2\cup K_3$, a contradiction. So, $N_{\{1,2\}}=\emptyset$. Similarly, $N_{\{3,4\}}=\emptyset$. Since $G$ is $(P_3\cup P_2)$-free, we have that $N_S=\emptyset$ if $|S|=1$. Therefore, $N=N_{\{1,4\}}\cup N_{\{1,5\}}\cup N_{\{2,4\}}\cup N_{\{2,5\}}\cup N_{\{1,2,4\}}\cup N_{\{1,2,5\}}\cup N_{\{1,4,5\}}\cup N_{\{2,4,5\}}\cup N_{\{1,2,4,5\}}$. Moreover, since $G$ is $(P_3\cup P_2, P_2\cup K_3)$-free,  we have that $G[M]$ is $(P_3, K_3)$-free, which implies that each component of $G[M]$ is a vertex or an edge.

If $N_{\{1,4\}}$ has two adjacent vertices $x_1$ and $x_2$, then $\{v_1,v_4,x_1,x_2,v_5\}$ induces a kite, a contradiction. So, $N_S$ is independent if $S\subseteq\{1,2,4,5\}$ and $|S|=2$. Similarly, $N_S$ is independent if $S\subseteq\{1,2,4,5\}$ and $|S|=3$. Let $J_2=N_{\{1,4\}}\cup N_{\{1,5\}}\cup N_{\{2,4\}}\cup N_{\{2,5\}}$ and $J_3=N_{\{1,2,4\}}\cup N_{\{1,2,5\}}\cup N_{\{1,4,5\}}\cup N_{\{2,4,5\}}$. We will prove that 

\begin{equation}\label{eqa-4}
	\mbox{$J_2$ is anticomplete to $J_3$}.
\end{equation}
Suppose not. Without loss of generality, there exists a vertex $y_1\in J_2$ such that $y_1\sim v_3$. If $y_1\in  N_{\{1,4\}}\cup N_{\{2,4\}}$, then $\{v_1,y_1,v_3,v_4,v_5\}$ or $\{v_2,y_1,v_3,v_4,v_5\}$ induces a kite, a contradiction. So, $y_1\in  N_{\{1,5\}}\cup N_{\{2,5\}}$. But now, $\{v_1,v_2,y_1,v_3,v_5\}$ induces a kite, a contradiction. So, $J_2$ is anticomplete to $J_3$. This proves (\ref{eqa-4}).

By (\ref{eqa-4}), we have that $\chi(G[N])\le\chi(G[J_2])+\chi(G[J_3])\le4$. If $v_3$ is not complete to $N_{\{1,2,4,5\}}$, let $y_2\in N_{\{1,2,4,5\}}$ such that $v_3\not\sim y_2$, then $\{y_2,v_1,v_2,v_3,v_5\}$ induces a kite, a contradiction. So, $V(Q)$ is complete to $N_{\{1,2,4,5\}}$. In particular, $\{v_1,v_2,v_3\}$ is complete to $N_{\{1,2,4,5\}}$, and thus $\omega(G[N_{\{1,2,4,5\}}])\le\omega(G)-3$.

Suppose $\omega(G)=3$. Then $N_{\{1,2,4,5\}}=\emptyset$, which implies that $\chi(G)\le\chi(G[N])+\chi(G[M\cup\{v_1,v_2,v_4,v_5\}])\le 4+2=6=2\omega(G)$, a contradiction. So, $\omega(G)\ge4$. Now, $\chi(G)\le\chi(G[N])+\chi(G[M\cup\{v_1,v_2,v_4,v_5\}])\le f(\omega(G)-3)+4+2\le2(\omega(G)-3)+6=2\omega(G)$, a contradiction. 

This proves Claim \ref{c2}.\qed

\begin{claim}\label{c3}
	$G$ is $(K_1\cup K_3)$-free.
\end{claim}
\pf Suppose not. Let $Q$ be an induced $K_1\cup K_3$ in $G$ with $V(Q)=\{u,v_1,v_2,v_3\}$ such that $v_1v_2v_3v_1$ is a triangle. By Lemma \ref{nonadjcent}, there exists a vertex $u'\in V(G)$ such that $u'\sim u$ and $u'\not\sim v_1$. By Claim \ref{c1} and \ref{c2}, $G$ is ($P_2\cup K_3$, hammer)-free.

If $u'$ is anticomplete to $\{v_2,v_3\}$, then $\{u,u',v_1,v_2,v_3\}$ induces a $P_2\cup K_3$, a contradiction. If $u'$ is complete to $\{v_2,v_3\}$, then $\{u,u',v_1,v_2,v_3\}$ induces a kite, a contradiction. So, $|N(u')\cap\{v_2,v_3\}|=1$. But now, $\{u,u',v_1,v_2,v_3\}$ induces a hammer, a contradiction. This proves Claim \ref{c3}.\qed

\medskip

\textbf{\noindent{\em Proof of Theorem~\ref{kite}} :} By Lemma \ref{C3}, we may assume that $\omega(G)\ge3$. Now, by Lemma \ref{K3} and Claim \ref{c3}, we have that $\chi(G)\le2\omega(G)$, a contradiction. This completes the proof of Theorem \ref{kite}.\qed

\medskip

Actually, by the proof above, we have the following proposition.

\begin{proposition}\label{p1}
	$\chi(F)\le2\omega(F)$ if $F$ is a $(P_3\cup P_2, K_1\cup K_3)$-free graph.
\end{proposition}

\section{($P_3\cup P_2$, hammer)-free graphs}
In this section, we consider ($P_3\cup P_2$, hammer)-free graphs. We may always assume that $G$ is a ($P_3\cup P_2$, hammer)-free graph such that $\chi(G')\le \omega^2(G')$ for every induced subgraph $G'$ of $G$ different from $G$, and $\chi(G)>\omega^2(G)$.  By Lemma \ref{C3}, we may assume that $\omega(G)\ge3$.

We will complete the proof of Theorem \ref{kite} by the following claim.

\begin{claim}\label{c3-1}
	$G$ is $(P_2\cup K_3)$-free.
\end{claim}
\pf Suppose not. Let $Q$ be an induced $P_2\cup K_3$ in $G$ with $V(Q)=\{u_1,u_2,v_1,v_2,v_3\}$ such that $v_1v_2v_3v_1$ is a triangle. We will prove that

\begin{equation}\label{eqa-3-1}
	\mbox{$N(u_1)=N(u_2)$}.
\end{equation}

Suppose not. Without loss of generality, let $u'\in N(u_1)\setminus N(u_2)$. If $|N(u')\setminus\{v_1,v_2,v_3\}|\ge2$, let $u'\not\sim v_1$ and $u'\not\sim v_2$, then $\{u_1,u_2,u',v_1,v_2\}$ induces a $P_3\cup P_2$, a contradiction. So, $|N(u')\cap\{v_1,v_2,v_3\}|\ge2$, let $u'\sim v_1$ and $u'\sim v_2$, then $\{u_1,u_2,u',v_1,v_2\}$ induces a hammer, a contradiction. Therefore, $N(u_1)=N(u_2)$. This proves (\ref{eqa-3-1}).

Let $N=N(\{u_1,u_2\})$ and $M=M(\{u_1,u_2\})$. By (\ref{eqa-3-1}), we have that $\{u_1,u_2\}$ is complete to $N$, which implies that $\omega(G[N])\le\omega(G)-2$. Since $G$ is $(P_3\cup P_2)$-free, it follows that $G[M]$ is $P_3$-free, and thus $\chi(G[M])\le\omega(G)$. Now, $\chi(G)\le \chi(G[N])+\chi(G[M\cup \{u_1,u_2\}])\le f(\omega(G)-2)+\omega(G)\le(\omega(G)-2)^2+\omega(G)\le\omega^2(G)$ as $\omega(G)\ge3$, a contradiction.

This proves Claim \ref{c3-1}.

\textbf{\noindent{\em Proof of Theorem~\ref{hammer}} :} By Claim \ref{c3-1}, we have that $G$ is $(P_2\cup K_3)$-free. Let $C=\{v_1,v_2,\dots,v_\omega\}$ be a maximum clique in $G$. We divide $V(G)\setminus N[v_1]$ as follows : $$A_2=\{v\in V(G)\setminus N[v_1] | v\not\sim v_2\},$$ $$A_i=\{v\in V(G)\setminus N[v_1] | v\not\sim v_i, v\not\in\cup_{j=2}^{i-1} A_j, 3\le i\le\omega\}$$

For $2\le i\le\omega$, since $A_i$ is anticomplete to $\{v_1,v_i\}$, we have that $G[A_i]$ is $P_3$-free. Consequently, $G[A_i]$ is $K_3$-free as $G$ is $(P_2\cup K_3)$-free, which implies that each component of $G[A_i]$ is a vertex or an edge. So, $\chi(G[\cup_{i=2}^{\omega} A_i])\le2(\omega(G)-1)$. 

Let $B=V(G)\setminus( \cup_{i=2}^{\omega} A_i\cup N[v_1])$. By the definition of $A_i$, we have that $B$ is complete to $C\setminus\{v_1\}$. So, $B$ is independent. Now, $\chi(G)\le \chi(G[N])+\chi(G[\cup_{i=2}^{\omega} A_i])+\chi(G[B])\le f(\omega(G)-1)+2(\omega(G)-1)+1\le(\omega(G)-1)^2+2\omega(G)-1=\omega^2(G)$, a contradiction. 

This proves Theorem \ref{hammer}.\qed

\medskip

In reality, by the proof above, we have the following proposition.

\begin{proposition}\label{p2}
	$\chi(F)\le\omega^2(F)$ if $F$ is a $(P_3\cup P_2, P_2\cup K_3)$-free graph.
\end{proposition}

\section{($P_3\cup P_2, C_5$)-free graphs}
In this section, we consider ($P_3\cup P_2, C_5$)-free graphs. Let $f(x)=\frac{3x^2+x}{2}$. We may always assume that $G$ is a ($P_3\cup P_2, C_5$)-free graph such that $\chi(G')\le f(\omega(G'))$ for every induced subgraph $G'$ of $G$ different from $G$, and $\chi(G)>f(\omega(G))$.  By Lemma \ref{C3}, we may assume that $\omega(G)\ge3$.

\medskip

\textbf{\noindent{\em Proof of Theorem~\ref{C5}} :} Let $v\in V(G)$, $\omega(G)=\omega$ and $A_i=\{u\in N(v) | \omega(G[M_{N^{\ge2}(v)}(u)])=i\}$, where $i=0,1,2,\dots,\omega$. Let $A'=\cup_{j=3}^\omega A_j$. We will prove that 

\begin{equation}\label{eqa-4-1}
	\mbox{$A'$ is a clique}.
\end{equation}

Suppose $A'$ has two nonadjacent vertices $v_1$ and $v_2$. Let $H_i$ be a component of $G[M_{N^{\ge2}(v)}(v_i)]$ with $\omega(H_i)=\omega(G[M_{N^{\ge2}(v)}(v_i)])$, for $i\in\{1,2\}$. By the definition of $A'$, we have that $\omega(H_1)\ge3$ and $\omega(H_2)\ge3$. For $i\in\{1,2\}$, let $Q_i$ be a triangle in $H_i$ with $V(Q_i)=\{x_i,y_i,z_i\}$.

If $|N(v_1)\setminus\{x_2,y_2,z_2\}|\ge2$, let $v_1\not\sim x_2$ and $v_1\not\sim y_2$, then $\{v_1,v,v_2,x_2,y_2\}$ induces a $P_3\cup P_2$, a contradiction. So, $|N(v_1)\cap\{x_2,y_2,z_2\}|\ge2$, and by symmetry, we may assume that $v_1\sim x_2$ and $v_1\sim y_2$. Similarly, we may suppose that $v_2\sim x_1$ and $v_2\sim y_1$. To forbid an induced $P_3\cup P_2$ on $\{v,v_1,x_2,x_1,y_1\}$, we have that $x_2\sim x_1$ or $x_2\sim y_1$. But now, $\{v,v_1,x_2,x_1,v_2\}$ or $\{v,v_1,x_2,y_1,v_2\}$ induces a $C_5$, a contradiction. So, $A'$ is a clique. This proves (\ref{eqa-4-1}).

Suppose $A_0\cup A_1\cup A_2=\emptyset$. Then $N(v)=A'$ and $|A'|\le\omega-1$. Let $v'\in A'$. Then $G[M(v')\cap N^{\ge2}(v)]$ is $P_3$-free as $G$ is $(P_3\cup P_2)$-free, and $\omega(G[N(v')\cap N^{\ge2}(v)])\le\omega-1$. So, $\chi(G)\le \chi(G[N(v)])+\chi(G[M(v')\cap N^{\ge2}(v)])+\chi(G[N(v')\cap N^{\ge2}(v)])\le \omega-1+\omega+f(\omega-1)=f(\omega-1)+2\omega-1\le f(\omega)$, a contradiction. 

Therefore, we suppose $A_0\cup A_1\cup A_2\ne\emptyset$. Let $C$ be a maximum clique in $G[A_0\cup A_1\cup A_2]$ with $C=\{t_1,t_2,\dots,t_{\omega_0}\}$. Let $I_i=C\cap A_i$ for $i\in\{0,1,2\}$. Let $D=\{u\in N^2(v) | u$ is complete to $C$$\}$. Let $B_0=\emptyset$ and $B_i=M_{N^2(v)\setminus D}(t_i)\setminus\cup_{j=0}^{i-1}B_j$.

For $i\in\{1,2,\dots,\omega_0\}$, if $i\in I_0$, then $B_i=\emptyset$ as $A_0$ is complete to $N^2(v)$. If $i\in I_1$, then $B_i$ is independent by the definition of $A_1$. That is to say, if $i\in I_0\cup I_2$, then $\chi(G[B_i])\le1$. Suppose $i\in I_3$. Since $G$ is $(P_3\cup P_2)$-free and $vt_i$ is an edge, it follows that $G[B_i]$ is $P_3$-free, and thus $G[B_i]$ is a union of cliques. By the definition of $A_2$, we have that $\chi(G[B_i])\le2$. Therefore, $\chi(G[N^2(v)\setminus D])\le2\omega_0$.

Since $G[N^{\ge3}(v)]$ is $P_3$-free, we have that $\chi(G[N^{\ge3}(v)])\le\omega$. By (\ref{eqa-4-1}), we see that $\chi(G[A'\cup N^{\ge3}(v)])\le\omega$. Now, $\chi(G)\le \chi(G[N(v)])+\chi(G[N^2(v)])+\chi(G[N^{\ge3}(v)])\le \chi(G[A_0\cup A_1\cup A_2])+\chi(G[A'])+\chi(G[D])+\chi(G[N^2(v)\setminus D])+\chi(G[N^{\ge3}(v)])\le f(\omega_0)+f(\omega-\omega_0)+2\omega_0+\omega\le f(\omega_0)+f(\omega-\omega_0)+3\omega-2\le f(1)+f(\omega-1)+3\omega-2\le f(\omega)$, a contradiction.

This proves Theorem \ref{C5}.\qed

\section{($P_3\cup P_2, K_4$)-free graphs}
In this section, we consider ($P_3\cup P_2, K_4$)-free graphs. By Lemma \ref{C3}, we may assume that $G$ contains a triangle. Let $G$ be a ($P_3\cup P_2, K_4$)-free graph, we will complete the proof of Theorem \ref{K4} by the two following claims.

\begin{claim}\label{c5-1}
	If $G$ contains an induced $2K_3$, then $\chi(G)\le6$.
\end{claim}
\pf Let $Q$ be an induced $2K_3$ in $G$ with $V(Q)=\{v_1,v_2,v_3,u_1,u_2,u_3\}$ such that $C=v_1v_2v_3v_1$ is a triangle. For a subset $S\subseteq\{1,2,3\}$, we define $N_S=\{v|v\in N(\{C\}),$ and $v_i\sim v$ if and only if $i\in S\}$. Note that $N(C)=N_{\{1\}}\cup N_{\{2\}}\cup N_{\{3\}}\cup N_{\{1,2\}}\cup N_{\{1,3\}}\cup N_{\{2,3\}}$ as $G$ is $K_4$-free.

Suppose $N_{\{1\}}\ne\emptyset$. Let $v'\in N_{\{1\}}$. If $|N(v')\setminus \{u_1,u_2,u_3\}|\ge2$, let $u_1\not\sim v'$ and $u_2\not\sim v'$, then $\{v_1,v_2,v',u_1,u_2\}$ induces a $P_3\cup P_2$, a contradiction. So, $|N(v')\cap \{u_1,u_2,u_3\}|\ge2$, which implies that $|N(v')\cap \{u_1,u_2,u_3\}|=2$ as $G$ is $K_4$-free. By symmetry, we may assume that $v'\sim u_1$, $v'\sim u_2$ and $v'\not\sim u_3$.  Now, $\{v',u_2,u_3,v_2,v_3\}$ induces a $P_3\cup P_2$, a contradiction.  Therefore, $N_{\{1\}}=\emptyset$, and by symmetry $N_{\{2\}}=N_{\{3\}}=\emptyset$.

Since $G$ is $K_4$-free, it follows that $N_{\{1,2\}}$, $N_{\{1,3\}}$ and $N_{\{2,3\}}$ are all independent. Moreover, $M(C)$ is $P_3$-free as $G$ is $(P_3\cup P_2)$-free, and thus $\chi(G[M(C)])\le3$ as $G$ is $K_4$-free. Therefore, $\chi(G)\le \chi(G[N(C)])+\chi(G[M(C)])\le3+3=6$. This proves Claim \ref{c5-1}.\qed

\begin{claim}\label{c5-2}
	If $G$ is $2K_3$-free and contains an induced $P_2\cup K_3$, then $\chi(G)\le6$.
\end{claim}
\pf Let $Q$ be an induced $P_2\cup K_3$ in $G$ with $V(Q)=\{v_1,v_2,v_3,u_1,u_2\}$ such that $C=v_1v_2v_3v_1$ is a triangle. For a subset $S\subseteq\{1,2,3\}$, we define $N_S=\{v|v\in N(\{C\}),$ and $v_i\sim v$ if and only if $i\in S\}$. Note that $N(C)=N_{\{1\}}\cup N_{\{2\}}\cup N_{\{3\}}\cup N_{\{1,2\}}\cup N_{\{1,3\}}\cup N_{\{2,3\}}$ as $G$ is $K_4$-free.

Let $v'\in N_{\{1\}}$. If $v'$ is anticomplete to $\{u_1,u_2\}$, then $\{v',v_1,v_2,u_1,u_2\}$ induces a $P_3\cup P_2$, a contradiction. If $v'$ is adjacent to exactly one element of $\{u_1,u_2\}$, say $u_1$, then $\{u_1,u_2,v',v_2,v_3\}$ induces a $P_3\cup P_2$, a contradiction. So, $v'$ is complete to $\{u_1,u_2\}$. Since $G$ is $K_4$-free, we have that $N_{\{1\}}$ is independent. Consequently, $N_{\{1\}}\cup N_{\{2\}}\cup N_{\{3\}}$ is independent.

Since $G$ is $K_4$-free, it follows that $N_{\{1,2\}}$, $N_{\{1,3\}}$ and $N_{\{2,3\}}$ are all independent. Moreover, $M(C)$ is $(P_3,K_3)$-free as $G$ is $(P_3\cup P_2, 2K_3)$-free, and thus $\chi(G[M(C)])\le2$. Therefore, $\chi(G)\le \chi(G[N(C)])+\chi(G[M(C)])\le1+3+2=6$. This proves Claim \ref{c5-2}. \qed

\medskip

\textbf{\noindent{\em Proof of Theorem~\ref{K4}} :} By Claim \ref{c5-1} and \ref{c5-2}, we may assume that $G$ is $(P_2\cup K_3)$-free. By Proposition \ref{p2}, we have that $\chi(G)\le\omega^2(G)\le9$. This proves Theorem \ref{K4}.\qed


\begin{thebibliography}{99}
\bibitem{BA18} A. P. Bharathi, S. A. Choudum, Colouring of ($P_3\cup P_2$)-free graphs, Graphs Comb., 34 (2018) 97-107.

\bibitem{BM1976} J.~A. Bondy, U.~S.~R. Murty, Graph Theory with Applications,
MacMillan, London, 1976.

\bibitem{BBS19} C. Brause, B. Randerath, I. Schiermeyer, and E. Vumar, On the chromatic number of $2K_2$-free
graphs, Disc. Appl. Math., 253 (2019) 14–24.

\bibitem{BMS19} C. Brause, M. Gei{\ss}er, I. Schiermeyer, Homogeneous sets, clique-separators, critical graphs, and optimal $\chi$-binding functions, Disc. Appl. Math., 320 (2022) 211-222.

\bibitem{CHM21} K. Cameron, S. Huang, and O. Merkel, An optimal $\chi$-bound for ($P_6$, diamond)-free graphs, J. of Graph Theory, 97 (2021) 451-465.


\bibitem{CK22} A. Char, T. Karthick, Optimal chromatic bound for ($P_3\cup P_2$, $\bar{P_3\cup P_2}$)-free graphs, arXiv preprint arXiv:2205.07447, 2022.


\bibitem{CK10} S. A. Choudum, T. Karthick, Maximal cliques in $(P_2\cup P_3, C_4)$-free graphs, Disc. Math., 310 (23) 3398-3403, 2010.

\bibitem{ER59}  P. Erd\"{o}s, Graph theory and probability, Can. J. Math., 11 (1959) 34-38.

\bibitem{CRSR06} M. Chudnovsky, N. Robertson, P. Seymour, R. Thomas, The strong perfect graph theorem, Ann. of Math., 164 (2006) 51-229.

\bibitem{CS10} M. Chudnovsky, P. Seymour, Claw-free graphs VI. Colouring, J. Comb. Theory. Ser. B, 100 (2010) 560-572.

\bibitem{GM22}  M. Gei{\ss}er,  Colourings of $P_5$-free graphs, PhD thesis, 2022.

\bibitem{Gy75} A. Gy\'{a}rf\'{a}s, On Ramsey covering-numbers, Infinite and Finite Sets (Colloq., Keszthely, 1973; dedicated to P. Erd\"{o}s on his 60th birthday), Vol. \uppercase\expandafter{\romannumeral2}, Colloq. Math. Soc. Janos Bolyai, Vol. 10, North-Holland, Amsterdam, pp. 801-816 (1975).

\bibitem{KM18} T. Karthick, S. Mishra, Chromatic bounds for some classes of $2K_2$-free graphs, Disc.
Math., 341 (2018) 3079–3088.

\bibitem{LLW23} R. Li, J. Li, and D. Wu, Optimal chromatic bound for ($P_3\cup P_2$, house)-free graphs, arXiv preprint arXiv:2308.05442, 2023.

\bibitem{LLWD23} R. Li, J. Li, and D. Wu, A tight linear chromatic bound for ($ P_3\cup P_2, W_4 $)-free graphs, arXiv preprint arXiv:2308.08768, 2023.

\bibitem{PA22} A. Prashant, P. Francis, S.F. Raj, $\chi$-binding functions for some classes of $(P_3\cup P_2) $-free graphs, arXiv preprint arXiv:2203.06423, 2022.

\bibitem{PA23} A. Prashant, S. Francis Raj, and M. Gokulnath, Bounds for the chromatic number of some
$pK_2$-free graphs, Disc. Appl. Math., 336 (2023) 99–108.

\bibitem{PAM23} A. Prashant, S. Francis Raj, and M. Gokulnath, Linear $\chi$-binding functions for ($P_3\cup P_2$, gem)-free graphs, arXiv preprint arXiv:2305.11757, 2023.

\bibitem{RS04} B. Randerath, I. Schiermeyer, Vertex colouring and forbidden subgraphs-a survey, Graphs Comb., 20 (2004) 1-40.

\bibitem{SR19} I. Schiermeyer, B. Randerath, Polynomial $\chi$-binding functions and forbidden induced subgraphs: a survey, Graphs Comb., 35 (2019) 1-31.

\bibitem{SS20} A. Scott, P. Seymour, A survey of $\chi$-boundedness, J. of Graph Theory, 95 (2020) 473-504.

\bibitem{Su81} D.P. Sumner, Subtrees of a graph and chromatic number, in: The Theory and Applications of Graphs, John Wiley \& Sons, New York, 1981, pp. 557-576.

\bibitem{Wa80} S. Wagon, A bound on the chromatic number of graphs without certain induced subgraphs,
J. Comb. Theory. Ser. B, 29 (1980) 345–346.

\bibitem{WZ22} X. Wang, D. Zhang, The $\chi$-Boundedness of $(P_3\cup P_2)$-Free Graphs, J. of Math., 2022. https://doi.org/10.1155/2022/2071887

\bibitem{WX23} D. Wu, B. Xu, Coloring of some crown-free graphs,  Graphs Comb., 39 (2023) 106.

\end{thebibliography}
\end{document}